\documentstyle[leqno,11pt]{article}

\newtheorem{Thm}{Theorem}[section]
\newtheorem{Lemma}[Thm]{Lemma}
\newtheorem{Prop}[Thm]{Proposition}
\newtheorem{Cor}[Thm]{Corollary}

\newcommand{\contr}{\,\rule{.1in}{.5pt}\rule{.5pt}{1.5mm}\,\,}
\newcommand{\surj}{\,\longrightarrow\hspace{-12pt}
                   \longrightarrow\,\,}
\newcommand{\downsurj}{\ \makebox[0in]
 {\raisebox{-2pt}{$\downarrow$}}\makebox[0in]{$\downarrow$}\ }

\newcommand{\llim}{\displaystyle \lim_
{\stackrel{\displaystyle\longleftarrow}{i}}}

\newcommand{\pf}{\ \\ \noindent {\bf Proof.\ \ }}
\newcommand{\qed}{\ $\displaystyle\Box$\\ \ \par}

\newcommand{\rk}{\refstepcounter{Thm}
\ \\ \noindent {\bf Remark \theThm }\ \ }

\newenvironment{Def}{\refstepcounter{Thm}
\ \\ \noindent {\bf Definition \theThm }\ \ }{\\ \ \par}

\newcommand{\II}{{\rm II}}
\newcommand{\boxtimes}{\makebox[0pt][l]{$\times$}
\raisebox{-1pt}{$\Box$}}
\newcommand{\scriptboxtimes}{\makebox[0pt][l]{$\scriptstyle\times$}
\makebox{$\scriptstyle\Box$}\,}

\renewcommand{\H}{{\cal H}}
\renewcommand{\d}{\mbox{\boldmath $d$}}
\newcommand{\C}{{\bf C}}
\newcommand{\D}{{\bf D}}
\newcommand{\R}{{\bf R}}
\newcommand{\T}{{\bf T}}
\newcommand{\Z}{{\bf Z}}
\newcommand{\HH}{{\bf H}}
\newcommand{\End}{\mbox{\rm End}}
\newcommand{\Aut}{\mbox{\rm Aut}}
\newcommand{\Hom}{\mbox{\rm Hom}}
\newcommand{\uu}{\overline{\cal U}^{(2)}}

\title{Infinite-dimensional Lie algebras\\
and the period map for curves}
\author{Yakov Karpishpan
\thanks{Supported in part by the NSF grant
DMS-9102233}}
\date{{\em March} 1994}

\begin{document}

\maketitle

\begin{abstract} We compute higher-order differentials of the period
map for curves and show how they factor through the corresponding
higher Kodaira-Spencer classes. Our approach is based on the
infinitesimal equivariance of the period map, due to Arbarello and De
Concini \cite{AD}.
\end{abstract}

\

\begin{flushright}
\em To the memory of Boris Moishezon
\end{flushright}

\

A holomorphic map between two complex manifolds
$\Phi:S\longrightarrow M$ induces, for every $n$, a map between
sheaves of differential operators of order $\leq n$ on $S$ and
$M$, ${\cal D}_S^{(n)}\longrightarrow{\cal D}_M^{(n)}$.
In Algebraic Geometry it is slightly more convenient to replace
${\cal D}_S^{(n)}$ with the sheaf of the $n^{th}$ order tangent
vectors $T^{(n)}S:={\cal D}_S^{(n)}/{\cal O}_S$. The $n^{th}$
differential of $\Phi$ is the corresponding map
$$
d^n\Phi:T^{(n)}S\longrightarrow T^{(n)}M\ .
$$

In general, it is hard to describe $T^{(n)}S$ or $T^{(n)}M$, and
little can be said about $d^n\Phi$. However, when $\Phi$ is the
period map arising from a family of complex algebraic varieties or
K\"{a}hler manifolds over $S$, and $M$ is an appropriate period
domain $\D$, it turns out that despite the transcendental nature
of $\Phi$, the differentials $d^n\Phi$ admit an algebraic
description.

In this paper we study higher differentials of the period map
$\Phi:S\longrightarrow\D$ associated to a miniversal deformation of
a complete curve $X$. Specializing to curves allows us to use constructions
and facts so far unavailable in other cases.
To simplify the exposition, we restrict our computations to $d^2\Phi$,
which already presents all features found in higher-order cases. This
also seems to be the most important case for potential applications.
The results for higher-order cases are summarized in Theorem
\ref{thm:main:higher}.

Our first main result is a description of the second
differential of $\Phi$ at $0\in S$,
$$
d_0^2\Phi:T^{(2)}_0S\longrightarrow T^{(2)}_{\Phi(0)}\D\ ,
$$
in terms intrinsic to $X$.

An essential difficulty one encounters is the lack of a nice
description for the second tangent space to the
period domain, $T^{(2)}_F\D$, comparable tothe well-known identification
$T_F\D=\Hom^{(s)}(F,H/F)$.
We bypass this problem by observing in Proposition
\ref{split} that there is a natural
splitting
$$
T^{(2)}_F\D=T_F\D\oplus S^2T_F\D\ ,
$$
and that $d_0^2\Phi$ splits accordingly:
$$
d_0^2\Phi=\ell\oplus \sigma\ .
$$
Here $\sigma$ is the {\em symbol part} of $d_0^2\Phi$, which is
simply $S^2d_0\Phi$ composed with the projection
$T_0^{(2)}S\rightarrow S^2T_0S$.

The object of our attention is the {\em linear part}of $d_0^2\Phi$,
$$
\ell:T_0^{(2)}S\longrightarrow T_{\Phi(0)}\D\ ,
$$
and in Theorem \ref{Thm:main} we explain how to compute it. The result may be
expressed in its simplest form as follows. Let $\zeta$ and $\xi$ be some
Kodaira-Spencer representatives for two tangent vectors to $S$ at 0.
Then$\zeta\otimes\xi$ represents,
 in a precise sense defined in the paper, a second-order tangent vector to $S$
at 0, and $\ell(\zeta\otimes\xi)$ is given by the map
$$
\omega\longmapsto\xi\contr\pounds_{\zeta}\omega
$$
regarded as an element of
$$
T_{\Phi(0)}\D=\Hom^{(s)}(H^0(X,\Omega^1_X),H^1(X,{\cal O}_X))\ .
$$

 It was shown in \cite{K1} that the same formulacomputes the {\em second
fundamental form} of  $\Phi$ (see \cite{CGGH}),
$$
\II: S^2T_0S\longrightarrow T_{\Phi(0)}\D/\,image\,(d_0\Phi)\ .
$$
In fact, $\II\equiv\ell\bmod image\,(d_0\Phi)$.
The very important refinement here is that $\II(\zeta\otimes\xi)$ is already
determined by the {\em symbol} of $\zeta\otimes\xi$ in $S^2T_0S$,whereas
$\ell(\zeta\otimes\xi)$ (and, hence, $d_0^2\Phi(\zeta\otimes\xi)$) involves the
full second-order tange
nt vector $[\zeta\otimes\xi]\in T_0^{(2)}S$.

The second main result of this paper is a cohomological interpretation of
$d_0^2\Phi$.
A well-known result of Griffiths \cite{Gri} states that the first differential
of the period map is given by cup product with the Kodaira-Spencer
class of the deformation.Recently there has appeared a series of
papers \cite{BG,EV,R1,R2,HS} renewing the study of higher-order
deformation theory and, in particular, introducing higher
Kodaira-Spencer classes $\kappa_n$. Pursuing an analogy with
Griffiths' result, we showed in \cite{K2} that the second
fundamental form of $\Phi$ depends only on $\kappa_2$ (more
precisely, on $\kappa_2$ modulo the image of $\kappa_1^2$).
Here this fact receives a new proof as a corollary of Theorem
\ref{Interpret:coho}, which states that $\ell$, and hence $d_0^2\Phi$, factors
through $\kappa_2$. In fact, Theorem \ref{Interpret:coho} gives more: it brings
$d_0^2\Phi$ in closer agreement wi
th $d_0^2\Phi$ by explicitly displayinga kind of a cup product computing $\ell$
on the cochain level.

In \cite{K1} and \cite{K2} the main technical
tool was Archimedean cohomology --- an infinite-dimensional
replacement for the Hodge structure of $X$. Here we use a different
infinite-dimensional object: the ``extended Hodge structure" of the
curve $X$, leading to the ``extended period map" $\widehat{\Phi}$.
Arbarello and De Concini introduced these notions in their paper
\cite{AD}, which serves as a point of departure for this work.

One starts with a curve $X$, a
point $p$ on $X$, and a formal parameter near $p$, $u:\widehat{\cal
O}_{X,p}\longrightarrow\C[[z]]$. We may say that both the usual and the
extended Hodge structures on $X$ encode information about regular
1-forms on $X-p$: one through their periods, the other through
their Laurent expansions at $p$.

The advantage of working with the extended period map is that it is
easier to bring in a basic fact of moduli theory for curves that
emerged in recent years --- that the moduli of curves are (locally)
infinitesimally uniformized by a very simple, though
infinite-dimensional, Lie algebra.We are referring to the {\em
Witt} Lie algebra of formal vector fields on a punctured disk,
$\d=\C((z))d/dz$,
whose central extension is the more famous Virasoro algebra.

The key observation of \cite{AD} used in this paper is that the
extended Hodge structures are also infinitesimally uniformized by
an infinite-dimensional Lie algebra, denoted ${\bf sp}(\H')$, and
that there is a Lie algebra homomorphism
$$
\varphi:\d\longrightarrow{\bf sp}(\H')
$$
making the extended, and hence the usual period map
{\em infinitesimally equivariant}. This means that the differential
of the period map $\Phi$ associated to a miniversal deformation of
$X$ over a sufficiently small base $S$ is induced by the (very
simple) Lie-algebraic object $\varphi$ (see diagram
(\ref{phi:induce:dPhi})). This sets the stage for computing higher
differentials of $\Phi$, which are evidently induced by filtered
pieces of the enveloping algebra morphism corresponding to $\varphi$.

The plan of the paper is as follows. Section 1 contains notation
and definitions pertaining to infinitesimal uniformization.
Section 2 is devoted to the proof of Theorem
\ref{thm:S:infhom} on infinitesimal uniformization of moduli of
curves. This result is due to\cite{BMS,Ko,BS}, cf. \cite{N,TUY}, but
the author was unable to find a complete proof in the literature.
Section 3 introduces extended Hodge structures and discusses
infinitesimal uniformization for the period domain. Section 4
reviews the construction of the extended Hodge structure
associated to a given curve and in doing so defines the extended
period map. In section 5 all of this comes
together in Theorem \ref{thm:equi} on infinitesimal equivariance
of the period map. The material in sections 3 through 5 is largerly
due to Arbarello and De Concini \cite{AD}, though we found it necessary to fill
in a number of details.

Sections 6 and 7 develop our main results. Proposition
\ref{split} exploits the usual (i.e. finite-dimensional)
uniformization of the period domain $\D$ under $Sp(2g,\C)$ to
obtain a splitting of the second tangent space to $\D$.
Then lemmas \ref{how:Sp:acts} and \ref{two:unifs} study the
infitesimal uniformization of $\D$ described in section 3 and
compare the two uniformizations on the second-order
level. As a corollary, we obtain in Theorem \ref{Thm:main}
a formula for the linear part of $d^2_0\Phi$, $\ell$.
The map $\ell$ also determines the second
fundamental form II of the VHS. We have computed
II in \cite{K1} in a different way, and in Theorem \ref{Thm:II} we
show that the two approaches agree.

Section 7 is devoted to showing that the second differential of
the period map, as well as its linear part $\ell$ introduced in
section 6, factors through the second Kodaira-Spencer class
$\kappa_2$. This requires an excursion into the recent description
of the second tangent space at $[X]$ to the moduli of curves as a
certain cohomology group on $X\times X$. In fact, in
Theorem \ref{Interpret:coho} we give what
amounts to a cohomological interpretation of the map $\ell$. In Corollary
\ref{II:factors} we obtain another proof that II factors through $\kappa_2$.
This was already proved
in \cite{K2} by a more general but less explicit method.

Section 8 provides higher-order analogues of the statements in
sections 6 and 7.

\newpage

\section{Notation and some preliminaries}
\label{sect:notat}

The results discussed in this paper owe their explicitness to a
very concrete object, the field of Laurent power series
$$
\H=\C((z))=\C[[z]][z^{-1}]\ .
$$
Most of the time we will regard it merely as an
infinite-dimensional vector space. It has several distinguished
subspaces:
$$
\H_+=\C[[z]]\ \ \ {\rm and} \ \ \
\H_-= \mbox{the span of negative powers of}\ z\ .
$$
Thus $\H=\H_+\oplus\H_-$ . Also,
$$
\H_+'=z\H_+, \ \ \ {\rm and} \ \ \ \H'=\H'_+\oplus \H_-\ .
$$

\begin{Def}
$<f,g>={\rm Res}_{z=0}fdg$\ .
\end{Def}
This is a symplectic form on $\H$, non-degenerate on $\H'$.

\begin{Def}
$${\bf sp}(\H')=\{\alpha\in
\End(\H')|<\alpha(x),y>+<\alpha(y),x>=0\ \ \mbox{\rm for all}\
x,y\in\H'\} \ .$$
\end{Def}

\noindent{\bf Facts:} (a) $\H$ (and, hence, $\H'$) is a
topological vector space with the $z$-adic topology.

(b) ${\bf sp}(\H')$ is isomorphic to the completion
$\widehat{S^2}(\H')$ of $S^2(\H')$ --- the symmetric square of
$\H'$, where $S^2(\H')$ embeds in ${\bf sp}(\H')$ via
$$
hk \mapsto \{x\mapsto <h,x>k+<k,x>h\}\ .
$$

Finally, $\d=\H\frac{d}{dz}$ will denote the {\em Witt} Lie
algebra of formal vector fields on a punctured disc. Its central
extension is the more famous Virasoro algebra.
We will also use $\d_+=\H_+\frac{d}{dz}$.

The Lie algebras $\d$ and ${\bf sp}(\H')$ will appear in the
setting of the following

\begin{Def}
\label{def:inf-hom}
A Lie algebra $L$ {\em acts by vector fields} on
a manifold $M$ if there is a homomorphism (or
anti-homomorphism) of Lie algebras
$$
L\longrightarrow \Gamma(M,\Theta_M)\ .
$$
If the composed map to the tangent space of $M$ at a point $x$
$$
L\longrightarrow \Gamma(M,\Theta_M)\longrightarrow T_xM
$$
is surjective for each $x\in M$, $M$ is called {\em
infinitesimally homogeneous}, and one says that $L$ provides an
{\em infinitesimal uniformization} for $M$.
\end{Def}

We will write $\Omega^1_X$ and $\omega_X$ interchangeably when $X$
is a curve. In turn, ``curve" will mean a complex algebraic curve.
We will consider the classical topology and the analytic structure
on $X$ only when dealing with the cohomology of $X$ with
coefficients in $\Z$ or $\C$, and when using the exponential
sequence in the proof of Lemma \ref{lemma:Lambda}, where we write
$X^{an}$.

We will also use the following notation: if ${\bf g}$ is a Lie
algebra, then ${\cal U}{\bf g}$ is its universal enveloping
algebra, and $\overline{\cal U}{\bf g}:={\cal U}{\bf g}/\C$.
${\cal U}^{(k)}{\bf g}$ (respectively, $\overline{\cal U}^{(k)}{\bf
g}$) will denote the elements of order$\leq k$ in the natural
filtration of ${\cal U}{\bf g}$ (respectively, $\overline{\cal
U}{\bf g}$).

Finally, we will use the fact that Lie algebras of endomorphisms of
Hodge structures, with or without polarization, carry a Hodge
structure of their own, always of weight 0. E.g. if $H=\oplus
H^{p,q}$ is a HS and ${\bf g}=\End(H)$, then the Hodge
decomposition of ${\bf g}$ is
$$
{\bf g}=\oplus{\bf g}^{-k,k}\ ,
$$
where ${\bf g}^{-k,k}=\{f\in\End(H)\ |\ f(H^{p,q})\subset
H^{p-k,q+k}\ \ \forall p,q\}$.

\section{Infinitesimal uniformization of moduli of curves}

Let $X$ be a complete curve of genus $g\geq 2$, $p$ --- a point
on $X$, and
$$
u:\widehat{\cal O}_{X,p}\stackrel{\cong}{\longrightarrow}\H_+
$$
--- a formal local coordinate at $p$; $u$ extends to an
isomorphism of fields of fractions, also defining the
obvious monomorphisms
$$
\Theta_{X,p}\hookrightarrow \H_+\frac{d}{dz}\subset\d\ , \ \ \
\Omega^1_{X,p}\hookrightarrow \H_+dz , \ \ \
\Omega_X^1(*p)_p\hookrightarrow \H dz , \ \ \ {\rm etc.},
$$
all of which will also be denoted by $u$. As any point on a
complete curve, $p$ is an ample divisor on $X$. Therefore,
$X-p$ is an affine open set in $X$. Choose an affine
neighborhood $V$ of $p$ in $X$. Then $\{(X-p), V\}$ is an affine
covering of $X$, suitable for computing the \v{C}ech cohomology
of $X$ with coefficients in a coherent sheaf. Thus we have an
exact sequence
\begin{eqnarray}
\label{seq:basic}
 & & \\ & &
0\rightarrow\Gamma(X-p,\Theta_X)\oplus\Gamma(V,\Theta_X)
\stackrel{\delta}{\rightarrow}\Gamma(V-p,\Theta_X)
\stackrel{\pi}{\rightarrow}H^1(X,\Theta_X)\rightarrow 0\ .
\nonumber
\end{eqnarray}
Exactness on the left is a consequence of $H^0(X,\Theta_X)=0$,
which, in turn, follows from the assumption $g\geq 2$.

The sheaf $\Theta_X$ is filtered by the subsheaves
$\Theta_X(-ip)$ of vector fields vanishing at $p$ to an order
$\geq i$ \ ($i\geq 0$). This induces a decreasing filtration
$P^i$ on spaces of sections over $V$, and hence over $X-p$ and
$V-p$. The \v{C}ech differential $\delta$ is strictly compatible
with $P^{\bullet}$, and so is the projection $\pi$, once
$H^1(X,\Theta_X)$ receives the induced filtration $P^{\bullet}$
from $\Gamma(V-p, \Theta_X)$. Therefore, the sequence
(\ref{seq:basic}) remains exact when reduced modulo $P^i$.

\begin{Lemma} The maps
$$
u_i:\Gamma(V-p, \Theta_X)/P^i\longrightarrow\d/z^i\d_+
$$
induced via the identification $u$ are isomorphisms for each
$i>0$.
\end{Lemma}
\pf Suppose the points $q_1,\ldots,q_m$ constitute the complement of $V$ in
$X$, and let $Q$ be the effective divisor $q_1+\ldots+q_m$. Since $Q$ is ample,
for $N$ sufficiently large $H^1(X,\Theta_X(NQ-ip))$ will vanish. We may also
assume $deg{\cal L}\ge
q g-1$. With this choice of $N$, set ${\cal L}=\Theta_X(NQ-ip)$. Then the
Riemann-Roch Theorem gives
$$
H^0(X,{\cal L}(kp))=\deg {\cal L}+k+1-g
$$
for each $k\geq 1$, which means that for each $k$ there exists a section of
$\Theta_V(-ip)$, regular on $V-p$ and with a pole of order exactly $k$ at $p$.
Thisimplies the surjectivity of $u_i$.

Now,
$$u:\Gamma(V-p,{\cal O}_X)\longrightarrow\H$$
is injective, since any regular function on $V-p$ is
completely determined by its Laurent expansion. And
$$u^{-1}(z^i\H_+)=\Gamma(V-p,{\cal O}_X(-ip))\ ,$$
implying that the $u_i$'s are injective too.

\qed

\begin{Cor}
$\llim\Gamma(V-p,\Theta_X)/P^i=\d$ .
\end{Cor}

\begin{Lemma}
\label{lemma:surj}
Passing to the inverse limit in the exact sequence obtained
from (\ref{seq:basic}) by reduction $\bmod P^i$ produces an exact
sequence
\begin{equation}
0\rightarrow
u(\Gamma(X,\Theta_X(*p)))\oplus\d_+\longrightarrow\d\longrightarrow
H^1(X,\Theta_X)\rightarrow 0
\label{seq:final}
\end{equation}
\end{Lemma}

\pf First we note that
$\Gamma(X-p,\Theta_X)=\Gamma(X,\Theta_X(*p))$ by \cite{Gro}. Also,
$P^i\Gamma(X,\Theta_X(*p))=0$ for all $i>0$, since $X$ supports no
non-zero global regular vector fields by virtue of the assumption
$g\geq 2$. Hence
$\Gamma(X,\Theta_X(*p))^{\wedge}=\Gamma(X,\Theta_X(*p))$.
 Second,
$$
\llim H^1(X,\Theta_X)/P^i=H^1(X,\Theta_X)\ ,
$$
because for all sufficiently large $i$
$$
H^1(X,\Theta_X)/P^i=H^1(X,\Theta_X)\ ;
$$
this is simply a consequence of $H^1(X,\Theta_X)$ being
finite-dimensional. Finally, inverse limits preserve the
exactness of (\ref{seq:basic}) $\bmod P^i$, because the directed
system
$$
\{A_i=(\Gamma(X-p,\Theta_X)\oplus\Gamma(V,\Theta_X))/P^i\}
$$
satisfies the Mittag-Leffler condition (see \cite{L}, III,
Prop. (9.3)):

{\em For each $n$, the decreasing sequence of images of natural
maps $\varphi_{mn}:A_m\rightarrow A_n\ \ (m\geq n)$ stabilizes.}

\noindent This is trivially so since all $\varphi_{mn}$ are
surjective in our situation. \qed

Assume now that $X$ moves in a flat family
\begin{eqnarray}
{\cal X} & \supset & X_t \nonumber\\
\pi\ \downarrow & & \downarrow \\
S & \ni & t\ \ ,\nonumber
\label{family:flat}
\end{eqnarray}
with a section ${\bf p}:S\rightarrow{\cal X}$ and a local
coordinate
$$
{\bf u}:\widehat{\cal O}_{{\cal X},{\bf p}}
\stackrel{\sim}{\longrightarrow}\Gamma(S,{\cal O}_S)
\otimes\H_+
$$
on $\cal X$ along $\bf p$, so that the restriction of $\bf u$
to $X_t$ provides a local formal coordinate $u_t$ near $p_t$.

 For each $t\in S$ one has an analogue of
(\ref{seq:final}). In particular, for each $t$ there is a
surjection$$
\d\longrightarrow H^1(X_t,\Theta_{X_t}) \ ;
$$
these glue together into a map

\begin{equation}
\d \longrightarrow\Gamma(S,R^1\pi_{*}\Theta_{{\cal X}/S}) \ .
\label{map:rel}
\end{equation}

Assume further that $S$ is a disc centered at $0$ in $\C^{3g-3}$,
and the family
\begin{equation}
\begin{array}{ccc}
{\cal X} & \supset & X \\
\pi\ \downarrow & & \downarrow \\
S & \ni & 0
\end{array}
\label{deform:X}
\end{equation}
is a miniversal deformation of the curve $X$. Then the
Kodaira-Spencer map of the family,
$$
\kappa:\Theta_S\longrightarrow R^1\pi_{*}\Theta_{{\cal X}/S}\ ,
$$
isan isomorphism. Composing its inverse with the map
in (\ref{map:rel}) yields a linear map
\begin{equation}
\lambda:\d\longrightarrow\Gamma(S,\Theta_S)\ .
\label{map:Lie}
\end{equation}

\begin{Lemma}
\label{lemma:act}
The map $\lambda$ in (\ref{map:Lie}) is an
anti-homomorphism of Lie algebras.
\end{Lemma}

\pf
$\cal X$ admits an acyclic Stein covering ${\cal W}=\{W_0,W_1\}$
with $W_0\cong S\times V$ and $W_1={\cal X}-{\bf p}\cong S\times
{X-{\bf p}}$. It follows from the proof of the previous lemma that
the map $\lambda$ fits in the commutative diagram
$$
\begin{array}{ccccc}
\Gamma(S,{\cal O}_S)\otimes\d & \hookleftarrow & \d &
\stackrel{\lambda}{\longrightarrow} & \Gamma(S,\Theta_S)\\ & & \\
j\uparrow\cong & &&  &
\cong\downarrow\kappa \\ & && & \\
\Gamma(W_0\cap W_1,\Theta_{{\cal X}/S})^{\wedge} & &
\stackrel{\textstyle c}{\surj} & & \Gamma(S,R^1\pi_{*}\Theta_{{\cal
X}/S}) \end{array}\ ,
$$
where $ ^{\wedge}$ indicates completion with respect to the
filtration by the order of vanishing along $\bf p$, and $j$ is the
isomorphism given by taking Laurent expansionsof relative vector
fields on $W_0\cap W_1$ along $\bf p$ via $\bf u$.

We begin by reviewing the
definition of $\kappa$. The Kodaira-Spencer map
$\kappa$ is the connecting morphism in the direct-image
sequence of the short exact sequence of ${\cal O}_{\cal X}$
modules
\begin{equation}
0\rightarrow\Theta_{{\cal X}/S}\longrightarrow\Theta_{\cal X}
\longrightarrow\pi^*\Theta_S\rightarrow 0\ \ .
\label{seq:KS}
\end{equation}
This contains an exact subsequence of $\pi^{-1}{\cal
O}_S$-modules
\begin{equation}
0\rightarrow\Theta_{{\cal X}/S}\longrightarrow
\widetilde{\Theta}_{\cal X}
\longrightarrow\pi^{-1}\Theta_S\rightarrow 0\ \ .
\label{subseq:KS}
\end{equation}
whose direct-image sequence also has $\kappa$ as a connecting
morphism (see \cite{BS} and also \cite{EV}). Furthermore,
(\ref{subseq:KS}) is an exact sequence of sheaves of Lie
algebras. The $\C$-linear brackets on $\widetilde{\Theta}_{\cal
X}$ and $\pi^{-1}\Theta_S$ are inherited from $\Theta_{\cal X}$
and $\Theta_S$, respectively. The bracket on $\Theta_{{\cal
X}/S}$ is even $\pi^{-1}{\cal O}_S$-linear.

We are ready to prove the lemma. Take any $\zeta,\xi\in\d$, and
let $Z=\lambda(\zeta),\ \Xi=\lambda(\xi)$. We wish to show that
$[\zeta,\xi]=-[Z,\Xi]$, where the first bracket istaken in
the Witt Lie algebra $\d$, and the second is in
$\Gamma(S,\Theta_S)$. The elements $\zeta$ and $\xi$ of $\d$,
which we identify with their pre-images under $j$, may be taken as
Kodaira-Spencer representatives of $Z$ and $\Xi$. Lift $Z$
to some sections of $\widetilde{\Theta}_{\cal X}$, $\zeta_0\in$ on
$W_0$ and $\zeta_1\in$ on $W_2$, and similarly for $\Xi$: $\xi_0\in
\Gamma(W_0,\widetilde{\Theta}_{\cal X})$, and $\xi_1\in
\Gamma(W_1,\widetilde{\Theta}_{\cal X})$. Then $\zeta_1-\zeta_0$
and $\xi_1-\xi_0$, with all terms restricted to $W_{01}=W_0\cap
W_1$, also give KS representatives for $Z$ and $\Xi$. In
particular,
$$
\zeta=\zeta_1-\zeta_0+\delta\theta\ ,
$$
and
$$
\xi=\xi_1-\xi_0+\delta\eta\ ,
$$
where $\theta$ and $\eta$ are some elements of $\check{C}^0({\cal
W},\Theta_{{\cal X}/S})^{\wedge}$. Then $[Z,\Xi]$ admits as its KS
representative the following expression, all terms of which are
restricted to $W_{01}$:
\begin{eqnarray}
\label{brackets}
\lefteqn{[\zeta_1,\xi_1]-[\zeta_0,\xi_0]=}\\
 & =& [\zeta_1,\xi_1]-
[\zeta_1-\zeta+\delta\theta,\xi_1-\xi+\delta\eta] \nonumber\\
 & =&
-[\zeta,\xi]+[\zeta_1,\xi]+[\zeta,\xi_1]+[\delta\theta,\xi_0]+
[\zeta_0,\delta\eta] \ .\nonumber
\end{eqnarray}
The Lie bracket of a section of
$\widetilde{\Theta}_{\cal X}$ with that of $\Theta_{{\cal X}/S}$
is again a section of $\Theta_{{\cal X}/S}$, which implies that
the last two terms in (\ref{brackets}) are in
$\delta\check{C}^0({\cal W},\Theta_{{\cal X}/S})$. We may assume
that $\bf u$ is induced by an isomorphism $u:\widehat{\cal O}_{X,p}
\rightarrow\H_+$ viathe identification $W_{01}\cong
S\times\{V-p\}$. The identification allows us to label some
vector fields on $W_{01}$ as horizontal or vertical. By
construction, $\zeta$ and $\xi$ are vertical and constant in the
horizontal direction. The fields $\zeta_1$ and $\xi_1$, on the
other hand, may be chosen to be horizontal and constant in the
vertical direction. Then $[\zeta_1,\xi]=[\zeta,\xi_1]=0$. Collecting
what is left of (\ref{brackets}), we conclude that $-[\zeta,\xi]$
is a Kodaira-Spencer representative for $[Z,\Xi]$,
which proves the lemma. \qed

Recallingdefinition \ref{def:inf-hom}, we may summarize lemmas
\ref{lemma:act} and \ref{lemma:surj} in the following theorem,
due to\cite{BMS,Ko,BS}, cf. \cite{N,TUY}.
\begin{Thm}
\label{thm:S:infhom}
For any curve $X$ of genus $g\geq 2$ the Witt Lie algebra $\d$
acts by vector fields on the base $S$ of a miniversal deformation
of$X$, making $S$ infinitesimally homogeneous.
\end{Thm}

\rk Theaction above clearly
depends on the choice of a point $p_t$ on each curve $X_t$, as
well as on a formal parameter $u_t$ at $p_t$. For our purposes
all these choices are equally good. More canonically, one may
consider the moduli space of triples $(X,p,u)$, encompassing
all possible choices of $p$ and $u$ on each $X$. The action of
$\d$ extends to such ``dressed" moduli spaces $\hat{\cal M}_g$,
making them also infinitesimally homogeneous. We will not need
these constructions, since the questions we study are local on
$\hat{\cal M}_g$.

\newpage

\section{Infinitesimal uniformization of period domains of weight
one}

By definition, a {\em Hodge structure of weight one} consists of a
lattice $\Lambda\cong\Z^{2g}$ and a decomposition of its
complexification $H=\Lambda\otimes\C$, $H=H^{1,0}\oplus H^{0,1}$,
such that $H^{1,0}=\overline{H^{0,1}}$. The {\em Hodge filtration}
$F^{\bullet}$ on $H$ is given by $F^0=H$, $F^1=H^{1,0}$,
$F^0=0$. The HS $(\Lambda,H,F^{\bullet})$ is {\em principally
polarized} if $\Lambda$ is equipped with a unimodular symplectic
form $Q(\ ,\ )$ such that $(u,v)=Q(\bar{u},v)$ is a
{\em positive-definite} Hermitian form on $H^{1,0}$ (and on
$H^{0,1}$).

The data $(\Lambda,H,F^{\bullet},Q)$ defines {\em a
principally-polarized abelian variety} $A=H^{0,1}/i(\Lambda)$,
where $i$ denotes the composition of the inclusion
$\Lambda\rightarrow\Lambda\otimes\C=H$ with the projection
$H=H^{1,0}\oplus H^{0,1}\rightarrow H^{0,1}$.

As is well-known, the space $\D$ of all Hodge structures
$(H,F^{\bullet})$ with a given lattice $\Lambda$ and polarization
$Q$ (={\em the period domain}) can be identified with the Siegel
upper half-space $\HH_g$ of complex symmetric $g\times g$ matrices
whose imaginary parts are positive-definite. The moduli space of
principally-polarized abelian varieties of dimension $g$, ${\cal
A}_g$, is a quotient of $\HH_g$ by the action of $Sp(2g,\Z)$.
Note that $\D$ is a homogeneous space for the group$Sp(2g,\R)$.

We wish to present $\D$ locally as an infinitesimally homogeneous
space for ${\bf sp}(\H')$.

\begin{Def}
An {\em extended Hodge structure} (of weight one) is a triple
$(Z,K_0,\Lambda)$, where $Z$ is a maximal isotropic subspace of
$\H'$ (with respect to the symplectic form $<\ ,\ >$), $K_0$ is a
codimension $g$ subspace of $Z$, and $\Lambda$ is a rank $2g$
lattice in $K^{\perp}_0/K_0$, subject to several conditions.
\end{Def}

First of all, $Z\cap\H'_+=0$. This implies the splittings
$\H'=Z\oplus\H'_+$ and
$$
H:=K^{\perp}_0/K_0=H^{1,0}\oplus H^{0,1}\ ,
$$
where $H^{0,1}=Z/K_0$, and $H^{1,0}=K^{\perp}_0\cap \H'_+$.

Let $Q$ be the bilinear form induced on $H$ by $\frac{1}{2\pi i}<\
 ,\ >$ on $\H'$. The
remaining conditions state that $H=\Lambda\otimes\C$, defining a
real structure on $H$, that $H^{1,0}=\overline{H^{0,1}}$ with
respect to this structure, and that $Q$ is unimodular on $\Lambda$.
Thus $(\Lambda, H, H^{1,0},H^{0,1},Q)$ is a principally-polarized
HS of weight one.

Arbarello and De Concini introduced an extended version of the
Siegel upper half-space, $\widehat{\HH}_g$, on which $Sp(2g,\Z)$
acts transitively, and the quotient manifold $\widehat{\cal A}_g$
parameterizes extended Hodge structures. The latter may also be
regarded as``extended abelian varieties" in view of the
following commutative diagram:
\begin{equation}
\begin{array}{ccc}
\widehat{\HH}_g & \longrightarrow & \widehat{\cal A}_g \\
 && \\
\downarrow & \swarrow & \downarrow \\
 && \\
\HH_g & \longrightarrow & {\cal A}_g
\end{array}
\label{diag:extended}
\end{equation}
The horizontal maps are quotients with respect to the
$Sp(2g,\Z)$-action. All spaces are manifolds (the top two are
infinite dimensional), except ${\cal A}_g$, which is a
$V$-manifold. Note that all maps in the upper triangle are smooth.

\begin{Prop}[\cite{AD}]
$\widehat{\cal A}_g$ is an infinitesimally
homogeneous space for ${\bf sp}(\H')$.
\label{prop:A:infhom}
\end{Prop}
Obviously, this also makes $\D=\HH_g$ {\em locally} infinitesimally
homogeneous for ${\bf sp}(\H')$. Let us work out the surjection
$$
{\bf sp}(\H')\longrightarrow T_H{\D}
$$
explicitly.
At any point $H=H^{1,0}\oplus H^{0,1}$ of $\D$,
$$
T_H\D=\Hom^{(s)}(H^{1,0}, H^{0,1})=S^2H^{0,1}\ .
$$
Suppose $H$ comes from an extended HS $(Z,K_0,\Lambda)$. Then any
$\alpha\in\End(\H')$ induces a map
\begin{equation}
H^{1,0}=K_0^{\perp}\cap \H'_+\longrightarrow\H' \ .
\label{sect:of:proj}
\end{equation}
We use the formulas $\H'=Z\oplus\H'_+$ and $K_0\cap \H'_+=0$ to observe
that
$$
H^{0,1}\cong Gr_F^{0} H= Z/K_0\cong \H'/K_0+\H'_+\ .
$$
Then (\ref{sect:of:proj}), composed with the natural projection
$$
\H'\longrightarrow\H'/K_0+\H'_+\ ,
$$
yields an element $a\in\Hom(H^{1,0},H^{0,1})$.

When $\alpha\in{\bf sp}(\H')$,
$$<\alpha(x),y>=-<x,\alpha(y)>\ ,$$
i.e. $<x,\alpha(y)>=<y,\alpha(x)>$ for all $x,y\in\H'$. Hence
$$
Q(x,a(y))=Q(y,a(x))
$$
for all $x,y\in H^{1,0}$, which means $a$ is {\em symmetric}:
$$
a\in\Hom^{(s)}(H^{1,0},H^{0,1})=S^2H^{0,1}\ .
$$
For reasons that will be clear later, we prefer $-a\in S^2H^{0,1}$.
Thus $\alpha\mapsto -a$ indeed defines a map
\begin{equation}
\label{res:rho}
\rho:{\bf sp}(\H')\longrightarrow T_H\D=
\Hom^{(s)}(H^{1,0},H^{0,1})=S^2H^{0,1}\ .
\end{equation}

\rk
\label{rho}
For further use we record that the above construction presents
the uniformizing map (\ref{res:rho}) as a restriction of a more
broadly defined map
\begin{equation}
\End(\H')\longrightarrow
\Hom(H^{1,0},H^{0,1})\ .
\end{equation}
Both maps will be denoted $\rho$.

\rk In view of (\ref{diag:extended}), Proposition
\ref{prop:A:infhom} implies that a sufficiently small open set $U$
in $\D$ is an infinitesimally homogeneous space under the action of
${\bf sp}(\H')$. However, the action is not unique --- it depends
on the choice of a lift from $U$ to $\widehat{\cal A}_g$.

\section{The extended period map}
Let $X$ be a complete smooth curve, $p$ --- a point on $X$, and
$u:\widehat{\cal O}_{X,p}\stackrel{\cong}{\longrightarrow}\H_+$
--- a formal local parameter at $p$. In this section we review
Arbarello and De Concini's construction associating an extended HS
$(Z,K_0,\Lambda)$ to the data $(X,p,u)$. When the triple $(X,p,u)$
varies in a flat family over some base $S$, this construction
defines ``an extended period map"
$$
\widehat{\Phi}: S\longrightarrow\widehat{\cal A}_g\ ,
$$
such that the usual period map $\Phi: S\longrightarrow \D$
naturally factors through $\widehat{\Phi}$.

\begin{Def}
$K_0:=u(\Gamma(X-p,{\cal O}_X))\cap\H'$\ .
\end{Def}

\noindent This is the same as putting $K_0=u(\Gamma(X-p,{\cal
O}_X))/\C$.

Note that $\Gamma(X-p,{\cal O}_X)=\Gamma(X,{\cal O}_X(*p))$ by a
theorem of Grothendieck \cite{Gro}, and that
$$
u:\Gamma(X-p,{\cal O}_X)\longrightarrow\H
$$
is injective. There are no non-constant regular functions on
$X$, and so $K_0\cap\H'_+=0$.

\begin{Lemma}
$H^1(X,{\cal O}_X)\cong \H/\H_++K_0$.
\label{H1:O}
\end{Lemma}

\pf
This follows from the exact sequence
$$
0\rightarrow\Gamma(V,{\cal O}_X)\oplus\Gamma(X-p,{\cal O}_X)\longrightarrow
\Gamma(V-p,{\cal O}_X)\longrightarrow H^1(X,{\cal O}_X)\rightarrow 0
$$
by completion with respect to the order-of-vanishing filtration $P^{\bullet}$
as in (\ref{lemma:surj}). \qed
Furthermore, $K_0$ is an isotropic
subspace of $\H'$, i.e. $K_0$ is contained in $K_0^{\perp}$, the
orthogonal complement of$K_0$ in $\H'$ with respect to the
symplectic form $<\ ,\ >$. We can be more specific about $K_0^{\perp}$.

\begin{Def}
$\Omega:=\{f\in\H'\ |\ df\in u(\Gamma(X-p,\Omega^1_X))\}$.
\end{Def}
We have $\Omega\cong\Gamma(X-p,\Omega^1_X)=\Gamma(X,\Omega^1_X(*p))$.

Now, Grothendieck's Algebraic De Rham Theorem \cite{Gro}, coupled wirh the
injectivity of the map $d:\H'\rightarrow\H dz$ and of $u$, gives
\begin{equation}
\Omega/K_0\cong
\frac{\Gamma(X,\Omega^1_X(*p))}{d\,\Gamma(X,{\cal O}_X(*p))} \cong
H^1(X-p,\C)=H^1(X,\C)\ .
\label{Groth:DR}
\end{equation}

\begin{Lemma}
$K_0^{\perp}=\Omega$.
\label{K:perp}
\end{Lemma}

\pf
If $f\in K_0$ and $g\in\Omega$, then $fdg$ is the Laurent expansion of a
globally defined one-form on $X$ with poles only at $p$. Then ${\rm Res}_0
fdg=0$, i.e. $<K_0,\Omega>=0$, and so $\Omega\subseteq K_0^{\perp}$. The
well-known duality theorem of Serr
e \cite{S} implies that the residue pairing induces a duality between
$H^0(X,\Omega_X^1)$ and $H^1(X,{\cal O}_X)$. The first of these groups is
isomorphic to $\Omega\cap\H'_+$, the second --- to
$$
\frac{\H}{\H_++K_0}=\frac{\Omega+\H}{\H_++K_0}\cong\frac{\Omega}{\Omega\cap(\H_++K_0)}= \frac{\Omega}{\Omega\cap\H_++K_0}\ .
$$
This implies that the residue pairing on $\Omega/K_0$ is non-degenerate.
Coupled with the earlier statements that $\Omega\subseteq K_0^{\perp}$ and
$K_0\subset\Omega^{\perp}$, we have $\Omega^{\perp}=K_0$ and
$(K_0^{\perp})^{\perp}\subseteq \Omega^{\perp}
$, which means that $(K_0^{\perp})^{\perp}=K_0$.

This, in turn, says that $<\ ,\ >$ is non-degenerate on $K_0^{\perp}/K_0$.
However, the pairing is 0 on $K_0^{\perp}\cap\H_+$ (since it is 0 on all of
$\H_+$), and on
$$
\frac{K_0^{\perp}}{(K_0^{\perp}\cap\H_+)+K_0}\cong
\frac{K_0^{\perp}+\H}{\H_++K_0}=\frac{\H}{\H_++K_0}\cong H^1(X,{\cal O}_X)\ .
$$
Then $K_0^{\perp}\cap\H_+$ must be dual to
$$
\frac{K_0^{\perp}}{(K_0^{\perp}\cap\H_+)+K_0}\cong H^1(X,{\cal O}_X)
$$
under the residue pairing on $K_0^{\perp}/K_0$, which implies
$K_0^{\perp}=\Omega$. \qed

\begin{Cor}
$K_0^{\perp}/K_0\cong H^1(X,\C)$.
\label{EHS:HS}
\end{Cor}

At this point we make the observation that the Laurent expansion
via $u$ at $p$ can be made well-defined not only for regular functions on a
punctured neighborhood of $p$, but also for sections of ${\cal
O}_X/\Z$:

\begin{Def}
 $K:=u(\Gamma(X-p,{\cal O}_X/\Z))\cap \H'$\ .
\end{Def}
Of course, by means of the exponential map, $\Gamma(X-p,{\cal
O}_X/\Z)$ may be regarded as a subspace of $\Gamma(X-p,{\cal
O}^*_{X^{an}})$.
In other words, $K$ consists of those $f\in\H'$ for which
$e^f$ lies in $u(\Gamma(X-p,{\cal O}^*_{X^{an}}))$. Obviously,
$K_0\subset K$.

Since the exterior derivative $d$ of a constant function is 0, $d$
is well-defined on ${\cal O}_X/\Z$, and (\ref{K:perp}) implies
that $K\subset K_0^{\perp}$.

\Def $\Lambda:=K/K_0$\ .

\begin{Lemma}
\label{lemma:Lambda}
The isomorphism (\ref{EHS:HS}):
$K_0^{\perp}/K_0\stackrel{\simeq}{\longrightarrow}H^1(X,\C)$ maps
$\Lambda$ onto $H^1(X,\Z)$; in particular,
$K_0^{\perp}/K_0\cong\Lambda\otimes\C$.
\end{Lemma}

\pf The starting point in identifying $H^1(X,\Z)$ with
$\Lambda$ is the exponential sequence (on $X^{an}$, of course)
\begin{equation}
\begin{array}{rcccccl}
 & & & & & {\cal O}^*_{X^{an}} & \\
 & & & & & e\,\uparrow\,\cong & \\
0\longrightarrow & \Z & \longrightarrow & {\cal O}_{X^{an}}
& \longrightarrow& {\cal O}_{X^{an}}/\Z & \longrightarrow 1\
\end{array}
\label{seq:exp}
\end{equation}
and its cohomology sequence
\begin{equation}\textstyle
\begin{array}{rcl}
H^1(X^{an},{\cal O}^*_{X^{an}})\  & & \Z\\
e\,\uparrow\,\cong\ \ \ \ \ \ \ \ \ \ & & \ \| \\
H^1(X^{an},\Z) \hookrightarrow
H^1(X^{an},{\cal O}_{X^{an}}) \rightarrow
H^1(X^{an}, {\cal O}_{X^{an}}/\Z) & \makebox[0pt]{$\,\surj$} &
H^2(X^{an},\Z)
\end{array}
\label{seq:exp:coho}
\end{equation}

But we also have an algebraic partial analogue of (\ref{seq:exp})
on $X$:
$$
0\longrightarrow\Z\longrightarrow{\cal O}_X\longrightarrow {\cal
O}_X/\Z\longrightarrow 1\ ,
$$
with the cohomology sequence
$$
\begin{array}{ccccccc}
0 & & & & & & 0 \\
\| & & & & & & \| \\
H^1(X,\Z) & \longrightarrow &
H^1(X,{\cal O}_X) & \longrightarrow & H^1(X,{\cal O}_X/\Z) &
\longrightarrow & H^2(X,\Z)
\end{array}
$$
mapping functorially to (\ref{seq:exp:coho}):
$$
\begin{array}{rcccl}
0\longrightarrow H^1(X^{an},\Z) \longrightarrow &
H^1(X^{an},{\cal O}_{X^{an}}) & \longrightarrow&
{\rm Pic}^0(X) & \longrightarrow 0\\
& \cong\ \uparrow & & \uparrow & \\
& H^1(X,{\cal O}_X) & \stackrel{\cong}{\longrightarrow} &
H^1(X,{\cal O}_X/\Z) & \ .
\end{array}
$$
The commutativity of the square implies that the right
verical arrow is surjective.

We also have the commutative ladder with exact columns
$$
\begin{array}{ccc}
0 & & \\
\uparrow & & \\
H^1(X,{\cal O}_X) & \stackrel{\cong}{\longrightarrow} &
H^1(X,{\cal O}_X/\Z)\\
\uparrow & & \uparrow \\
\Gamma(V-p,{\cal O}_X) & \longrightarrow &
\Gamma(V-p,{\cal O}_X/\Z)\\
\uparrow & & \uparrow \\
\Gamma(V,{\cal O}_X)\oplus\Gamma(X-p,{\cal O}_X) & \longrightarrow
& \Gamma(V,{\cal O}_X/\Z)\oplus\Gamma(X-p,{\cal O}_X/\Z) \\
\uparrow & & \uparrow \\
0 & & 0
\end{array}
$$
Again we note that the upper right vertical arrow must be
surjective.

Splicing the two diagrams, and completing with respect to the
order-of-vanishing filtration $P^{\bullet}$ as in (\ref{lemma:surj}),
we get
\begin{equation}
\begin{array}{rcccl}
 & 0 & & 0 & \\
& \uparrow & & \uparrow & \\
0\longrightarrow H^1(X^{an},\Z) \longrightarrow &
H^1(X^{an},{\cal O}_{X^{an}}) & \longrightarrow&
{\rm Pic}^0(X) & \longrightarrow 0\\
& \uparrow & & \uparrow & \\
& \H & = & \H & \\
& \uparrow & & \uparrow & \\
& \H_++K_0 & \longrightarrow & \H_++K & \\
& \uparrow & & \uparrow & \\
& 0 & & 0 &
\end{array}
\label{ladder:ABC}
\end{equation}

The lemma now follows by simple homological algebra. Consider the
vertical ladder in the above diagram as a short exact sequence of
three complexes
$$
0\longrightarrow A^{\bullet}\longrightarrow B^{\bullet}
\longrightarrow C^{\bullet}\longrightarrow 0\ .
$$
Then $H^0(C^{\bullet})=H^1(X^{an},\Z)$, $H^1(A^{\bullet})=K/K_0$,
and the connecting map in the corresponding
cohomology sequence is precisely the sought-after isomorphism
\begin{equation}
H^1(X^{an},\Z)\stackrel{\cong}{\longrightarrow}\Lambda=K/K_0\ .
\label{sought-after}
\end{equation}

It remains to show that this isomorphism is induced by theone in
(\ref{EHS:HS}). The map (\ref{sought-after}) factors through the monomorphism
$$
H^1(X^{an},\Z)\longrightarrow H^1(X^{an},{\cal O}_{X^{an},})\ ,
$$
which, in turn, factors through $H^1(X^{an},\C)$. And the vertical sequences in
(\ref{ladder:ABC}) may be amended as in the proof of Lemma \ref{K:perp}. Then
we arrive at the following variant of (\ref{ladder:ABC}):
$$
\begin{array}{rccl}\textstyle
 & H^1(X^{an},\C) & 0 & \ \ \ \ \ 0 \\
 & \nearrow\ \ \uparrow\ \ \searrow & \uparrow &\ \ \ \ \ \uparrow  \\
H^1(X^{an},\Z) & \longrightarrow &
H^1(X^{an},{\cal O}_{X^{an}}) & \rightarrow
{\rm Pic}^0(X)\rightarrow 0\\
& | & \uparrow & \ \ \ \ \ \uparrow  \\
& \ \ \ K_0^{\perp}\ = & K_0^{\perp} & =  \ K_0^{\perp}  \\
& \uparrow & \uparrow & \ \ \ \ \ \uparrow  \\
& \ \ \ K_0\ \rightarrow & (K_0^{\perp} \cap \H_+)+K_0 & \rightarrow
(K_0^{\perp} \cap \H_+)+K  \rightarrow\Lambda\\
& \uparrow & \uparrow & \ \ \ \ \ \uparrow \\
& 0 & 0 & \ \ \ \ \ 0
\end{array}
$$
With this diagram it is easy to trace the map $H^1(X^{an},\Z) \longrightarrow
\Lambda$ and see that it fits in the commutative square
$$
\begin{array}{ccc}
H^1(X^{an},\C) & \longrightarrow & K_0^{\perp}/K_0\\
\uparrow & & \uparrow\\
H^1(X^{an},\Z) & \longrightarrow & \Lambda
\end{array}
$$
with natural inclusions as the vertical arrows. \qed

\begin{Prop}
The isomorphism in (\ref{EHS:HS}) is {\em symplectic}, identifying
${\frac{1}{2\pi i}<\ ,\ >}$ on $\Lambda$ with the intersection form $Q(\ ,\ )$
on $H^1(X,\Z)$.
\end{Prop}

\pf The above proposition is established in \cite{AD}, following \cite{SW}, by
reasoning
similar to that in the proof of the Riemann reciprocity laws. Alternatively, we
can identify the residue pairing with the cup product
$$
H^0(X,\Omega_X)\otimes H^1(X,{\cal O}_X)\longrightarrow H^1(X,\Omega_X^1)\cong
H^2(X,\C)\cong \C\ ,
$$
as Serre suggests in \cite{S}, and then relate the cup product to the
intersection pairing. \qed

We now complete the identifications above to include the Hodge
structure. First, $U:=K_0^{\perp}\cap\H'_+$ is easily seen to be
mapped onto
$$
H^0(X,\Omega_X^1)=H^{1,0}(X)=F^1H^1(X,\C)
$$
by the isomorphism (\ref{EHS:HS}). Let $\overline{U}$ be the
complex conjugate of $U$ with respect to the real structure which
$\Lambda$ defines on $K_0^{\perp}/K_0=\Lambda\otimes\C$. Then
(\ref{EHS:HS}) identifies $\overline{U}$ with $H^{0,1}(X)$.
Finally, let $Z\subset K_0^{\perp}$ to be the pre-image of
$\overline{U}$ with respect to the projection
$$
K_0^{\perp}\surj K_0^{\perp}/K_0\ .
$$
It is easy to see that $Z\cap\H'_+=0$,
 $\H'=Z\oplus\H'_+$, and that $Z$ is a maximal isotropic subspace
of $\H'$.

To summarize, we have constructed an extended HS $(Z,K_0,\Lambda)$
out of the data $(X,p,u)$.

\section{Infinitesimal equivariance of the period map}
We will work with a miniversal deformation $\pi:{\cal
X}\rightarrow S$ of a complete smooth curve $X$ of genus $g\geq
2$, as in (\ref{deform:X}), with a sufficiently small contractible
open Stein manifold $S$ as its base.

It was shown in Theorem \ref{thm:S:infhom} that $S$ is an
infinitesimally homogeneous space for $\d$. Consider the usual
and the extended period maps on $S$:
\begin{eqnarray}
S & \stackrel{\widehat{\Phi}}{\longrightarrow} & \widehat{\cal
A}_g\nonumber\\
 & &\nonumber\\
\Phi\downarrow & j\nearrow & \downsurj \\
 & &\nonumber\\
U & \subset & \D=\HH_g\ .\nonumber
\label{maps:period}
\end{eqnarray}
Let $U$ be a neighborhood of $\Phi(0)$ in $\D$ containing the
image of $S$; we assume that $U$ is small enough to admit lifts to
$\widehat{\cal A}_g$. Choose the lift $j$ making the diagram
commutative (i.e. $j\circ\Phi=\widehat{\Phi}$ on $S$). This makes
$U$ an infinitesimally homogeneous space for ${\bf sp}(\H')$.

\Def $\varphi: \d \longrightarrow {\bf sp}(\H')$ is the
Lie-algebra homomorphism given by
$$
f\frac{d}{dz}\longmapsto\{g\mapsto fg'\ \ \ \forall g\in \H'\}\ .
$$

Using the identification ${\bf sp}(\H')\cong\widehat{S^2}(\H')$
(see Section \ref{sect:notat}), the map $\varphi$ may also be
written as
\begin{eqnarray*}
\varphi:\d & \longrightarrow & \widehat{S^2}(\H')\\
z^{k+1}\frac{d}{dz} & \longmapsto &
\frac{1}{2}\sum_{j\in\Z-\{0\}} z^{-j}z^{j+k}\ .
\end{eqnarray*}
We note that $\varphi$ is an irreducible representation of the
Witt algebra on $\H'$, described in \cite{KR}, (1.2), where it is
denoted $V'_{0,0}$.

The following is an adaptation of a theorem of Arbarello and De
Concini \cite{AD}.
\begin{Thm}
\label{thm:equi}
The period map $\Phi: S\longrightarrow U\subset\D$
is infinitesimally equivariant, i.e. there exists a commutative
diagram
\begin{equation}
\label{phi:induce:dPhi}
\begin{array}{ccc}
\d & \stackrel{\varphi}{\longrightarrow} & {\bf sp}(\H') \\
\downarrow & & \downarrow \\
\Gamma(S,\Theta_S) & \stackrel{d\Phi}{\longrightarrow} &
\Gamma(U,\Theta_{\D}) \ .
\end{array}
\end{equation}
The vertical arrows are Lie algebra
anti-homomorphisms, while the horizontal ones are Lie algebra
homomorphisms. The vertical arrows induce surjections onto $T_tS$
(respectively, $T_{H}\D$) for any point $t\in S$ (respectively,
$H\in U\subset\D$). \end{Thm}

\rk The vertical arrows are not unique.

\newpage

\section{The second differential of the period map}
\label{second:diff}
We continue with a miniversal deformation(\ref{deform:X}) of
$X$. Theorem \ref{thm:equi} allows one to calculate the various
differentials of the period map. We
begin by specializingdiagram (\ref{phi:induce:dPhi}) to $0\in S$:
\begin{equation}
\label{diag:equi:at_0}
\begin{array}{ccc}
\d & \stackrel{\varphi}{\longrightarrow} & {\bf sp}(\H') \\
\downsurj & & \downsurj \\
T_0S & \stackrel{d_0\Phi}{\longrightarrow} & T_{\Phi(0)}\D
\end{array}
\end{equation}
A well-known theorem of Griffiths \cite{Gri} factors $d_0\Phi$ as
\begin{equation}
\label{diag:Griffiths}
\begin{array}{ccl}
T_0S & \stackrel{d_0\Phi}{\longrightarrow} & T_{\Phi(0)}\D \\
\kappa\downarrow\cong & & \ || \\
H^1(\Theta_X) & \stackrel{\nu}{\longrightarrow} &
\Hom^{(s)}(H^0(\omega_X),H^1({\cal O}_X))\ ,
\end{array}
\end{equation}
where $\kappa$ is the Kodaira-Spencer isomorphism, and
$\nu$ is the map defined by the cup-product pairing
\begin{equation}
\label{cup}
H^1(\Theta_X)\otimes
H^0(\omega_X)\stackrel{\smile}{\longrightarrow} H^1({\cal O}_X)\ ,
\end{equation}
itself induced by the contraction pairing of sheaves
$\Theta_X\otimes\omega_X\stackrel{\contr}{\rightarrow}{\cal O}_X$.

Splicing (\ref{diag:equi:at_0}) and (\ref{diag:Griffiths}) yields
\begin{equation}
\label{diag:post-Griffiths}
\begin{array}{ccl}
\d & \stackrel{\varphi}{\longrightarrow} & {\bf sp}(\H')\\
\downsurj & & \ \downsurj{\scriptstyle \rho} \\
H^1(\Theta_X) & \stackrel{\nu}{\longrightarrow} &
\Hom^{(s)}(H^0(\omega_X),H^1({\cal O}_X))\ ,
\end{array}
\end{equation}
which we want to work out explicitly. As before, let $V$ be an
affine open set in $X$ containing $p$, so that $X-p$ and $V$ form
an affine covering of $X$. Let $\xi$ be a vector field on $V-p$
with $u(\xi)=f(z)\frac{d}{dz}$ in $\d$. Let $\omega$ be a global
holomorphic 1-form on $X$ with $u(\omega|_{V-p})=dg$ for some
$g\in\H'$. Then the cup-product pairing (\ref{cup}) gives
$$
[\xi]\smile[\omega]=[-\xi\contr(\omega|_{V-p})]\in H^1({\cal
O}_X)\ .
$$
We observe that the minus sign is built into $\rho$ (see
(\ref{rho})), and that
$$
u(\xi\contr(\omega|_{V-p}))=f\frac{d}{dz}\contr
dg=fg'=\varphi(f\frac{d}{dz})g\ ,
$$
which is how (\ref{diag:post-Griffiths}) and, indeed, the theorem
of Arbarello and De Concini (\ref{thm:equi}) is proved.

We would like to work out an equally explicit realization of the
second differential of $\Phi$ (higher-order cases are
similar). Our starting point is again Theorem \ref{thm:equi}. We
simply pass from Lie algebras to their (reduced) enveloping
algebras to obtain a commutative diagram
\begin{equation}
\label{diag:U:2}
\begin{array}{lcl}
\uu\d & \stackrel{\varphi^{(2)}}{\longrightarrow} & \uu{\bf sp}(\H')
\\ \downarrow & & \downarrow \\
\Gamma(S,\Theta^{(2)}_S) & \stackrel{d^2\Phi}{\longrightarrow} &
\Gamma(U,\Theta^{(2)}_{\D})\ ,
\end{array}
\end{equation}
where $\Theta^{(2)}={\cal D}^{(2)}/{\cal O}$ stands for the
second-order tangent sheaf, and $\uu$ is the notation introduced
in Section \ref{sect:notat}.

Again, to be precise, the maps emanating from the upper-left
corner reverse the order of products, while the remaining maps are
the second-degree parts of filtered ring homomorphisms.

Restricting to $0\in S$, we obtain
\begin{equation}
\label{equi:at_0:two}
\begin{array}{ccl}
\uu\d & \stackrel{\varphi^{(2)}}{\longrightarrow} &
\uu{\bf sp}(\H')\\
\lambda^{(2)}\downsurj\& & \ \downsurj\rho^{(2)} \\
T_0^{(2)}S & \stackrel{d^2_0\Phi}{\longrightarrow} &
T_0^{(2)}\D\ .
\end{array}
\end{equation}

\begin{Prop}
\label{split}
The second tangent space of the period domain $\D$ at the
point corresponding to a HS $(H,F^{\bullet})$ admits a canonical
splitting $$
T_F^{(2)}\D=T_F\D\oplus S^2T_F\D\ .
$$
\end{Prop}

\pf
Let ${\bf g}=\End(H)$ ($={\bf gl}(2g,\C)$), and ${\bf
s}={\bf sp}(H)$ (symplectic with respect to the polarization on
$H$). Then $\D$ is infinitesimally homogeneous under the action of
$\bf s$, and
$$
T_F\D\cong{\bf s}^{-1,1}\cong\Hom^{(s)}(F^1,H/F^1)\ .
$$
We also have a natural surjection $\uu{\bf s}\surj T_F^{(2)}\D$.
Its restriction to $\uu{\bf s}^{-1,1}$ is an isomorphism by reason
of dimension. But ${\bf s}^{-1,1}$ is an abelian Lie algebra, i.e.
$$
\uu{\bf s}^{-1,1}={\bf s}^{-1,1}\oplus S^2{\bf s}^{-1,1}\ .$$
\qed

In view
of Proposition \ref{split},
$d^2_0\Phi:T_0^{(2)}S\rightarrow T_{\Phi(0)}^{(2)}\D$ breaks up
into a direct sum of two components:

{\em the symbol map}
$$
\begin{array}{ccl}
T_0^{(2)}S & \stackrel{\sigma}{\longrightarrow} & S^2{\bf
s}^{-1,1}\\
(\Upsilon+\sum_iZ_i\Xi_i)|_0 & \longmapsto &
d_0\Phi(Z|_0)\otimes d_0\Phi(\Xi|_0)\ \ \ (\mbox{\rm order does
not matter here}),
\end{array}
$$
where $\Upsilon,Z_i,\Xi_i\in\Gamma(S,\Theta_S)$, and

{\em the linear part}
$$
\ell: T_0^{(2)}S\longrightarrow{\bf s}^{-1,1}\ .
$$
It is the linear part that is really interesting. A typical
second-order tangent vector to $S$ at 0,
$(\Upsilon+\sum_iZ_i\Xi_i)|_0$, is sent by $\ell$ to
$$
d_0\Phi(\Upsilon)+\sum_i\ell((Z_i\Xi_i)|_0)\ .
$$
Thus, it suffices to understand $\ell((Z\Xi)|_0)$ for
$Z,\Xi\in\Gamma(S,\Theta_S)$.

By surjectivity of $\lambda^{(2)}$ in
(\ref{equi:at_0:two}), we may assume that the vector fields $Z$
and $\Xi$ on $S$ are the images, respectively, of
some $f_1\frac{d}{dz}$ and $f_2\frac{d}{dz}$ in $\d$, under the map
$\d\rightarrow\Gamma(S,\Theta_S)$, whose restriction$\lambda$ is.
Then (\ref{equi:at_0:two}) implies that
\begin{equation}
\label{vw}
d_0^2\Phi((Z\Xi)|_0)=\rho^{(2)}\circ\varphi^{(2)}(f_1\frac{d}{dz}
f_2\frac{d}{dz})=\rho^{(2)}(vw) \,
\end{equation}
where $v=\varphi(f_1\frac{d}{dz})$ and
$w=\varphi(f_2\frac{d}{dz})$.

Now, the map$\rho^{(2)}$ in (\ref{equi:at_0:two})
is not induced by $\rho:{\bf sp}(\H')\rightarrow{\bf s}^{-1,1}$,
which was a restriction of the map, also denoted $\rho$ in
(\ref{rho}),
$$
\End(\H')\longrightarrow{\bf g}^{-1,1}\ .
$$
In fact, the maps $\rho$ are not even Lie algebra morphisms.
Nevertheless, there is a way to reduce$\rho^{(2)}$ to
$\rho$. This will require a more detailed understanding of the
infinitesimal action of ${\bf sp}(\H')$; in fact, we need to work
out how the group $Sp(\H')\subset\Aut(\H')$ acts on a neighborhood
of a point in $\widehat{\cal A}_g$.

Since any element in the group $\Aut(\H')$ may be written as
$I+\alpha$, where $\alpha\in\End(\H')$, we have the following map
from $\Aut(\H')$ to $\Aut(H)$:
$$
A\longmapsto I-\rho(\alpha), \ {\rm where}\ I+\alpha=A^{-1}\ .
$$
This map will be denoted $R$. So
\begin{equation}
R(A)=I-\rho(A^{-1}-I)\ .
\label{R}
\end{equation}

\noindent {\bf Caution:} $R$ is not a group homomorphism.

Let $(H,F^{\bullet}_t)$ be the Hodge structure corresponding to a
point in $U$ near $\Phi(0)$. The HS $(H,F^{\bullet}_t)$ comes from
an extended HS $(Z_t,K_{0,t},\Lambda_t)$. The assumption that $U$ is
small and infinitesimally homogeneous under ${\bf sp}(\H')$ implies
that there exists $A_t\in Sp(\H')$ such that $K_{0,t}$ and $Z_t$
are images under $A_t$ of $K_0$ and $Z$, respectively (we refer to
the components of the extended HS corresponding to $\Phi(0)$).

\begin{Lemma}
\label{how:Sp:acts}
In this situation $F_t^1=R(A_t)F^1$.
\end{Lemma}
\pf
We wish to compare the Hodge structures
$$
(H=K_0^{\perp}/K_0,F^1=K_0^{\perp}\cap\H'_+)
$$
and
$$
(H_t=K_{0,t}^{\perp}/K_{0,t},F_t^1=K_{0,t}^{\perp}\cap\H'_+)\ .
$$
To do so, we identify $H_t$ with $H$ by $A_t^{-1}$. Then the
comparison involves two subspaces of $H$:
$F^1=K_0^{\perp}\cap\H'_+$ and
$$
A_t^{-1}F_t^1=A_t^{-1}(K_{0,t}^{\perp})\cap A_t^{-1}(\H'_+)=
K_0^{\perp}\cap A_t^{-1}(\H'_+)\ .
$$

We regard $U$ as a subset of the Grassmannian
$$
Grass(F^1,H)=\Aut(H)/\{A\,|\,A(F^1)\subseteq F^1\}\ .
$$
Any element of $\Aut(H)$ may be written as $I+T$ for some
$T\in{\bf g}$, and if$I+T\in\{A\,|\,A(F^1)\subseteq F^1\}$,
then $T\in{\bf g}^{0,0}$. If some $I+T\in\Aut(H)$ moves
$F^1$ to $A_t^{-1}F^1_t$, then so does $I+T^{-1,1}$, where the
subscript refers to the $(-1,1)$-component of $T$ under the direct
sum decomposition ${\bf g}={\bf g}^{-1,1}\oplus{\bf g}^{0,0}\oplus{\bf
g}^{1,-1}$.
Thus we only need to find the map
$$
T^{-1,1}:H^{1,0}=\H'_+\cap K_0^{\perp}\longrightarrow
\H'/\H'_++K_0\cong H^{0,1}
$$
which measures deviation of $A_t^{-1}F^1_t$ from $F^1$. The
above formulas for $F^1$ and $A_t^{-1}F^1_t$ show that $T^{-1,1}$
is induced by
$$
A_t^{-1}:\H'_+\longrightarrow \H'/\H'_+\ .
$$
But if $A_t^{-1}=I+\alpha$, then
$$
\alpha:\H'_+\longrightarrow \H'/\H'_+
$$
induces the same $T^{-1,1}$. Recalling the definition of $\rho$
(\ref{rho}), this says that $T^{-1,1}=-\rho(\alpha)$. It remains
to consult the definition of $R$ (\ref{R}) and to abuse notation
by putting $F_t^1=A_t^{-1}F_t^1$. \qed

We are now able to establish the principal formula relating the
two infinitesimal uniformizations of $U$ on the second-order level.

\begin{Lemma}
\label{two:unifs}
Let $v,w\in\End(\H')$. Then
$$
\rho^{(2)}(vw)=\rho(v)\rho(w)-\rho(w\circ v)\ ,
$$
where $w\circ v$ denotes the composition law in $\End(\H')$.
\end{Lemma}
\pf
Let $V=\rho(v)$, $W=\rho(w)$. $V$ and $W$ are vectors in
$T_{\Phi(0)}\D\cong {\bf s}^{-1,1}$. For any $Z\in{\bf s}^{-1,1}$
we will write $\widetilde{Z}$ to denote the vector field on $U$
correspondingto $Z$ under the Lie algebra homomorphism
$$
{\bf s}^{-1,1}\longrightarrow\Gamma(U,\Theta_{\D})\ .
$$
In particular, $\widetilde{V}|_{\Phi(0)}=V$,
$\widetilde{W}|_{\Phi(0)}=W$. Let $f$ be any smooth function on
$U$. Then
\begin{eqnarray*}
\lefteqn{\rho^{(2)}(vw)f=} & & \\
 & = & \frac{\partial^2}{\partial t\partial s}|_0
      f\{R(\exp tv\circ\exp sw)\Phi(0)\}\\
 & = & \frac{\partial^2}{\partial t\partial s}|_0
      f\{[I-\rho(\exp(-sw)\circ\exp(-tv) -I)]\Phi(0)\}\\
& = & \frac{\partial^2}{\partial t\partial s}|_0
      f\{[I-\rho(-sw-tv+tsw\circ v+\ldots)]\Phi(0)\}\\
 & = & \frac{\partial^2}{\partial t\partial s}|_0
      f\{(I+tV+sW-ts\rho(w\circ v)+\ldots)\Phi(0)\}\\
 & = & \frac{d}{dt}|_0\left\{\frac{\partial}{\partial s}|_0
      f\{[(I+tV+o(t))+s(W-t\rho(w\circ v)+o(t))+o(s)]\Phi(0)\}
       \right\}\\
 & = & \frac{d}{dt}|_0\left\{
      [(\widetilde{W}-t(\rho(w\circ v))^{\sim}+o(t))f]
     \{(I+tV+o(t))\Phi(0)\}\right\}\\
 & = & \frac{d}{dt}|_0\biggl\{
      [\widetilde{W}f]\{(I+tV+\ldots)\Phi(0)\}-
     t[(\rho(w\circ v))^{\sim}f]\{(I+tV+\ldots)\Phi(0)\}+\\
 & & +o(t)   \biggr\}\\
 & = & V\widetilde{W}f-\rho(w\circ v)f\ .
\end{eqnarray*}\qed

Observe that $\rho$ takes
its values in ${\bf g}^{-1,1}$, which is an abelian Lie algebra.
Thus the lemma gives a splitting of
$$
\rho^{(2)}:\uu\End(\H')\longrightarrow\uu{\bf g}^{-1,1}={\bf
g}^{-1,1}
\oplus S^2{\bf g}^{-1,1}\ :
$$
$V\widetilde{W}=\rho(v)\rho(w)$ ispurely quadratic (=the symbol
part), and $-\rho(w\circ v)$ is the linear part. Going back to
(\ref{vw}),this implies that $\ell((Z\Xi)|_0)=-\rho(w\circ v)$,
which proves the following
\begin{Thm}
\label{Thm:main}
If $Z,\Xi\in \Gamma(S,\Theta_S)$ lift to $f_1\frac{d}{dz},
f_2\frac{d}{dz}\in\d$, then
the linear part of $d^2_0\Phi$,
$$
\ell:T_0^{(2)}S\longrightarrow T_{\Phi(0)}\D={\bf s}^{-1,1}\ ,
$$
sends
$(Z\Xi)|_0$ to the negative of the image under
$\rho:\End(\H')\rightarrow{\bf g}^{-1,1}$ of the composition in
$\End(\H')$ of $\varphi(f_1\frac{d}{dz})$ and
$\varphi(f_2\frac{d}{dz})$, in reverse order:
\begin{equation}
g\longmapsto f_2f'_1g'+f_1f_2g''\ .
\label{prescription}
\end{equation}
\end{Thm}

\rk A composition (in $\End(\H')$) of two elements of ${\bf
sp}(\H')$ need not be in ${\bf sp}(\H')$. In particular, it is not
a priori obvious that the image of $\ell$ is in
$$
\Hom^{(s)}(H^0(\omega_X),H^1({\cal
O}_X))\cong{\bf s}^{-1,1}\ .
$$

There is a better-known object which carries part of the
information contained in the linear part $\ell$ of the period
map's second differential. It is {\em the second fundamental form
of the VHS} of \cite{CGGH}, the map
$$
\II: T_0^{(2)}S/T_0S=S^2T_0S\longrightarrow
T_{\Phi(0)}\D/im\,(d_0\Phi)
$$
induced by $\ell$.

\begin{Thm}
\label{Thm:II}
The prescription (\ref{prescription}) for computing $\ell$ gives
a formula for \II, which coincides with that in \cite{K1}, $\S 6$:
\begin{equation}
\label{II}
\II(Z\otimes\Xi)=
\{\omega\mapsto\xi\contr\pounds_\zeta\omega\}\bmod im\,(d_0\Phi)
\end{equation}
for any$Z,\Xi\in T_0S$ with KS representatives
$\zeta,\xi\in\Gamma(V-p,\Theta_X)$, and $\omega\in H^0(X,\omega_X)$.
\end{Thm}

\pf Recall that a choice of a point $p$ on the curve $X$ and a local
parameter $z$ near $p$ allows one to represent
$\omega\in H^0(X,\omega_X)$ by some $g\in\H'$ with $dg=\omega$ near
$p$. The vectors $Z$ and $\Xi$ are the images under $\rho$ of
some $f_1\frac{d}{dz}$ and $f_2\frac{d}{dz}$ in $\d$,
i.e. $f_1\frac{d}{dz}$ and $f_2\frac{d}{dz}$ are the
Laurent expansions at $p$ of $\zeta$ and $\xi$, respectively.
Working out (\ref{II}) in terms of $z$, using the formula for the
Lie derivative
$$ \pounds_\zeta\omega=d\zeta\contr\omega+\zeta\contr
d\omega\ ,
$$
easily yields (\ref{prescription}) ($\bmod\
im\,(d_0\Phi)$),as was already done, in fact, in  \cite{K1}, $\S
6$. \qed
\rk
Thus the prescription for $\II$ given in \cite{K1} turns out to be
well-defined for second-order differential operators on $S$ and
not just their symbols, and the values given by that prescription
are not merely equivalence classes modulo $im\,(d_0\Phi)$.

\newpage

\section{Relation with the second Kodaira-Spencer class}
\label{rel:withKS2}
As explained at the beginning of the previous section, the first
differential of the period map is given bycup product with the
(first) Kodaira-Spencer class $\kappa=\kappa_1$ of the
deformation. In \cite{K2} we have shown that $\II$ depends only on
the {\em second} Kodaira-Spencer class $\kappa_2$ (more precisely,
on $\kappa_2\ \bmod\ im\,(\kappa_1)$) introduced recently in
\cite{BG}, \cite{EV} and \cite{R1}. In this section we will
explain, in the case of curves, how the full second differential
$$
d^2_0\Phi: T_0^{(2)}S\longrightarrow T_{\Phi(0)}^{(2)}\D
$$
factors through the second KS mapping
$$
\kappa_2: T_0^{(2)}S\longrightarrow \T_X^{(2)}\ .
$$

Let us recall first the construction of $\T_X^{(2)}$, the space of
second-order deformations of $X$. Our reference is \cite{R1} or
\cite{K2}.

Let $X_2$ denote the symmetric product of the curve $X$ with
itself;write
$$
g: X\times X\longrightarrow X_2
$$
for the obvious projection map, and $i:X\hookrightarrow X_2$ for the
inclusion of the diagonal. Then $\T_X^{(2)}=\HH^1(X_2,{\cal
K}^{\bullet})$, where ${\cal K}^{\bullet}$ is the sheaf complex on
$X_2$
$$
\begin{array}{ccc}
{\scriptstyle -1} & & {\scriptstyle 0}\\
(g_*(\Theta_X^{\scriptboxtimes 2}))^- & \stackrel{[\ ,\
]}{\longrightarrow} & i_*\Theta_X\ .
\end{array}
$$
Here $\boxtimes$ stands for the exterior tensor product on
$X\times X$, $(\ \ )^-$ denotes anti-invariants of the
$\Z/2\Z$-action, and the differential is the restriction to the
diagonal followed by the Lie bracket of vector fields.

In practice it seems easier to do the following. Letting
$C^{\bullet}$ denote the \v{C}ech cochain complex
$\check{C}^{\bullet}({\cal U},\Theta_X)$ of $\Theta_X$ with
respect to an affine covering $\cal U$ of $X$, one may compute
$\T_X^{(2)}$ as the cohomology of the simple complex associated to
the double complex
$$
\begin{array}{cccc}
{\scriptstyle 2} & (C^1\otimes C^1)^{(s)} & & \\
 & 1\otimes\delta\uparrow -\delta\otimes 1 & & \\
{\scriptstyle 1} & (C^0\otimes C^1 + C^1\otimes C^0)^-  &
\stackrel{[\ ,\ ]}{\longrightarrow} & C^1 \\
 & 1\otimes\delta\uparrow +\delta\otimes 1 & & \uparrow\delta \\
{\scriptstyle 0} & (C^0\otimes C^0)^- & \stackrel{[\ ,\
]}{\longrightarrow} & C^0 \\
& {\scriptstyle -1} & & {\scriptstyle 0}
\end{array}
$$
The superscripts $^{(s)}$ and $^-$ denote the invariants and the
anti-invariants, respectively, of the $\Z/2$-action.

Working with ${\cal U}=\{X-p,V\}$ and using Laurent expansions at
$p$, we may follow the proof of Lemma \ref{lemma:surj} and replace
$C^1$ with $\d$ and $C^0$ with the completion of its image in
$\d\oplus\d$. The resulting bicomplex still computes $\T_X^{(2)}$.
In particular, $\T_X^{(2)}$ is a quotient of
$\d\oplus(\d\otimes\d)^{(s)}$.

\begin{Lemma}
\label{ks2:rep}
Assume the vector fields $Z$ and $\Xi$ on $S$ are the images
of $\zeta,\xi\in\d$ under the infinitesimal uniformization map
$\lambda: \d\longrightarrow\Gamma(S,\Theta_S)$. Then
$$
\frac{1}{2}([\xi,\zeta]+\zeta\otimes\xi+\xi\otimes\zeta)\in
\d\oplus(\d\otimes\d)^{(s)}
$$
is a representative for $\kappa_2((Z\Xi)|_0)\in \T_X^{(2)}$.
\end{Lemma}

\pf The Kodaira-Spencer maps are compatible with the symbol map
$\T_X^{(2)}\longrightarrow S^2\T_X^1$ in the sense that there is a
commutative diagram
$$
\begin{array}{ccc}
T_0^{(2)}S & \stackrel{\kappa_2}{\longrightarrow} & \T_X^{(2)}\\
\downarrow && \downarrow\\
S^2T_0S & \stackrel{\kappa_1^2}{\longrightarrow} & S^2\T_X^1\ .
\end{array}
$$
Thus it is natural to look for a representative of
$\kappa_2((Z\Xi)|_0)$ of the form
$$
\theta+\frac{1}{2}(\zeta\otimes\xi+\xi\otimes\zeta)\in
\d\oplus(\d\otimes\d)^{(s)}
$$
for some $\theta\in \d$. The construction of $\kappa_2$ as the
connecting morphism in a certain long exact sequence
(\cite{EV,R1,R2}), presented more explicitly in \cite{K2}, offers
the following way to determine $\theta$. Working with a covering
${\cal W}=\{W_0,W_1\}$ of $\cal X$ as in the proof of Lemma
\ref{lemma:act}, and using the subsheaf $\widetilde{\Theta}_{\cal
X}$ of $\Theta_{\cal X}$ introduced in (\ref{subseq:KS}), let
$\zeta_0,\zeta_1$ be lifts of $Z$ to $\Gamma(W_0,
\widetilde{\Theta}_{\cal X})$, $\Gamma(W_1,
\widetilde{\Theta}_{\cal X})$. Write $\widetilde{\zeta}$ to denote
$\zeta_0+\zeta_1$ viewed as a cochain in $\check{C}^0({\cal W},
\widetilde{\Theta}_{\cal X})$. A slight modification of the proof
of Prop. 2 in \cite{K2} shows that $\theta$ should be cohomologous
(in $\check{C}^1({\cal W},\Theta_{{\cal X}/S})$) to
$$
\frac{1}{2}([\widetilde{\zeta},\xi]+[\xi,\widetilde{\zeta}])=
\frac{1}{2}([\zeta_0,\xi]+[\xi,\zeta_1])=
\frac{1}{2}[\xi,\zeta_1-\zeta_0]\ .
$$
But $\zeta_1-\zeta_0$ is cohomologous to $\zeta$. Hence we can take
$\theta=\frac{1}{2}[\xi,\zeta]$. \qed

In (\ref{split}) we explained how $T_{\Phi(0)}^{(2)}\D$
splits into $T_{\Phi(0)}\D\oplus S^2T_{\Phi(0)}\D$, with the
second differential of the period map breaking up accordingly:
$$
d_0^2\Phi=\ell\oplus \sigma\ .
$$

The symbol part factors through the square of the first KS class:
$$
\begin{array}{ccrcl}
T_0^{(2)}S && \stackrel{\sigma}{\longrightarrow} & &
S^2T_{\Phi(0)}\D \\
 & & & & \\
\downsurj & & {\scriptstyle (d_0\Phi)^2}\nearrow & & \ \ \uparrow
\nu_1^2\\
 & & & & \\
S^2T_0S & & \stackrel{\kappa_1^2}{\longrightarrow} & &
S^2\T_X=S^2\Hom^{(s)}(H^0(\omega_X),H^1({\cal O}_X))\ .
\end{array}
$$
This diagram may be directly obtained from (\ref{diag:Griffiths})
and carries no additional information.

Now to the linear part $\ell$ of $d_0^2\Phi$. ``Recall" the
canonical bijection $b$ given by the composition of the obvious
maps
$$
b: \d\oplus(\d\otimes\d)^{(s)}\hookrightarrow
\d\oplus(\d\otimes\d) \surj \uu\d\ .
$$

\begin{Lemma}
\label{lemma:canon.bij}
The canonical bijection $b$ fits in the commutative
square
$$
\begin{array}{ccc}
\uu\d & \stackrel{b}{\longleftarrow} &
\d\oplus(\d\otimes\d)^{(s)}\\
\downsurj & & \downsurj \\
T_0^{(2)}S & \stackrel{\kappa_2}{\longrightarrow} & \T_X^{(2)}
\end{array}
$$
with bijective horizontal arrows, and surjective vertical ones.
\end{Lemma}

\pf It suffices to show that if $\zeta,\xi\in\d$ lift the vector
fields $Z$ and $\Xi$ on $S$, then $b^{-1}(\zeta\xi)$ lifts
$\kappa_2((\Xi Z)|_0)$ under the projection
$$
\d\oplus(\d\otimes\d)^{(s)}\surj \T_X^{(2)}\ .
$$
In other words, we
must verify that $b^{-1}(\zeta\xi)$ is a KS representative for
$(\Xi Z)|_0\in T_0^{(2)}S$. From the definition of $b$ it easily
follows that
$$
b^{-1}(\zeta\xi)=
\frac{1}{2}([\xi,\zeta]+\zeta\otimes\xi+\xi\otimes\zeta)\ .
$$
This, together with Lemma \ref{ks2:rep}, implies our statement. \qed

\begin{Def}
We define $\nu_2:\T_X^{(2)}\longrightarrow T_{\Phi(0)}\D=
\Hom^{(s)}(H^0(\omega_X),H^1({\cal O}_X))$ as the composition
$$
\T_X^{(2)}\stackrel{\kappa_2^{-1}}{\longrightarrow}T_0^{(2)}S
\stackrel{\ell}{\longrightarrow}T_{\Phi(0)}\D\ .
$$
\end{Def}
Thus we have a commutative triangle
\begin{equation}
\begin{array}{ccccc}
T_0^{(2)}S & & \stackrel{\ell}{\longrightarrow} & &
T_{\Phi(0)}\D\\
 & & & & \\
 & \searrow\kappa_2 & & \nu_2\nearrow & \\
 & & & & \\
 & & \T_X^{(2)} & & \ .
\end{array}
\end{equation}

\begin{Thm}
\label{Interpret:coho}
 $\nu_2:\T_X^{(2)}=\HH^1({\cal
K}^{\bullet})\longrightarrow\Hom(H^0(\omega_X),H^1({\cal O}_X))$
is induced by the pairing
\begin{equation}
\HH^1({\cal K}^{\bullet})\otimes H^0(\omega_X)\longrightarrow
H^1({\cal O}_X)\ ,
\label{pairing:coho}
\end{equation}
defined on the \v{C}ech cochain level by the coupling
\begin{equation}
\label{pairing:cochain}
\begin{array}{cclcl}
(\check{C}^1(\Theta_X))^{\otimes 2} \oplus
\check{C}^1(\Theta_X) & \times\ & \check{C}^0(\omega_X) &
\longrightarrow & \check{C}^1({\cal O}_X) \\
(\zeta\otimes\xi+\upsilon) &
\times & \omega & \longmapsto &
\xi\contr\pounds_{\zeta}\omega-\upsilon\contr\omega\ .
\end{array}
\end{equation}
\end{Thm}

\pf It suffices to study the effect of
$\nu_2=\ell\circ\kappa_2^{-1}$ on an element
$x$ of $\T_X^{(2)}$ represented by
$$
\frac{1}{2}(\zeta\otimes\xi+\xi\otimes\zeta)+\upsilon\in
(\d\otimes\d)^{(s)}\oplus \d\ .
$$
According to Lemma \ref{ks2:rep},
$$
\kappa_2^{-1}(x)=(Z\Xi+\Upsilon-\frac{1}{2}[Z,\Xi])|_0\ ,
$$
where $Z,\Xi$ and $\Upsilon$ are the images of $\zeta,\xi$ and
$\upsilon$, respectively, under the uniformization map $\lambda:
\d\longrightarrow\Gamma(S,\Theta_S)$. Note that
$\lambda([\zeta,\xi])=-[Z,\Xi]$.

As explained in the previous section,
$\ell((Z\Xi+\Upsilon-\frac{1}{2}[Z,\Xi])|_0)$ is a map
$H^0(\omega_X)\longrightarrow H^1({\cal O}_X)$ given by
\begin{equation}
\omega\longmapsto\xi\contr\pounds_{\zeta}\omega
-(\upsilon-\frac{1}{2}[\zeta,\xi])\contr\omega\ .
\label{pre-Cartan}
\end{equation}
However, Cartan's identity gives
$$
\pounds_{\zeta}(\xi\contr\omega)-\xi\contr\pounds_{\zeta}\omega=
[\zeta,\xi]\contr\omega\ ,
$$
and on a curve $\pounds_{\zeta}(\xi\contr\omega)=
\zeta\contr\pounds_{\xi}\omega$. Hence the right-hand side of
(\ref{pre-Cartan}) equals
$$
\frac{1}{2}(\xi\contr\pounds_{\zeta}\omega+
\zeta\contr\pounds_{\xi}\omega)-\upsilon\contr\omega\ .
$$
\qed

\begin{Cor}
\label{II:factors}
The second fundamental form of the VHS, \II, factors through
$S^2\T_X$.
\end{Cor}

\noindent This fact was already proved in complete generality
(for a deformation of any compact K\"{a}hler manifold) in
\cite{K2}, using Archimedean cohomology.

We conclude with a diagram summarizing the relationships between
some of the maps discussed in this section:
\begin{equation}\begin{array}{ccccc}
 & & \ell & & \\
 & & & & \\
T_0^{(2)}S & \stackrel{\kappa_2}{\longrightarrow} &
\T^{(2)}_X & \stackrel{\nu_2}{\longrightarrow} &
T_{\Phi(0)}\D \\
\downsurj && \downsurj & & \downsurj\\
S^2T_0S & \stackrel{\kappa_1^2}{\longrightarrow} &
S^2\T_X &
\stackrel{\nu_2/\,im\,\nu_1}{\longrightarrow} & T_{\Phi(0)}\D/\,im\,\nu_1 \\
 & & & & \\
& & \II & &
\end{array}
\label{relationships}
\end{equation}

\section{The higher-order case}
We have the following analogues of the results in sections
\ref{second:diff} and \ref{rel:withKS2}. The proofs, which are
notationally cumbersome transcriptions of the $n=2$ case, are omitted.

\begin{Prop}
The $n^{th}$ tangent space of the period domain $\D$ at a point
corresponding to a HS $(H,F^{\bullet})$ admits a canonical
splitting
$$
T_F^{(n)}\D=T_F\D\oplus S^2T_F\D \oplus\ldots\oplus S^nT_F\D\ .
$$
\end{Prop}

\noindent The $n^{th}$ differential of the period map splits
accordingly:
$$
d_0^n\Phi=\ell_1^{(n)}+\ldots+\ell_n^{(n)}\ .
$$
E.g. what we called $\ell$ and $\sigma$ earlier are $\ell_1^{(2)}$
and $\ell_2^{(2)}$, respectively.

Thus it suffices to describe the $k^{th}$ component of $d_0^n\Phi$,
$$
\ell_k^{(n)}: T_0^{(n)}\D\longrightarrow S^k T_{\Phi(0)}\D\ .
$$

\begin{Thm}
If $Z_1,\ldots,Z_n\in\Gamma(S,\Theta_S)$ lift to
$\zeta_1,\ldots,\zeta_n\in\d$, then
\newcounter{bean}
\begin{list}{\rm\alph{bean})}{\usecounter{bean}}
\item $\ell_1^{(n)}$ sends $(Z_1\ldots Z_n)|_0$ to $(-1)^{n-1}$ times the image
under
$\rho:\End(\H')\longrightarrow{\bf g}^{-1,1}$ of the composition
in $\End(\H')$ of
$\varphi(\zeta_1),\ldots,\varphi(\zeta_n)$ in reverse
order;
\item $\ell_k^{(n)}$ is the sum, over all partitions of $k$, of the
symmetrized tensor products
$$
\overline{\bigotimes}_{\sum_i p_i=k}\ell_1^{(p_i)}\ ;
$$
\item $d_0^{n}\Phi$, as well as each $\ell_k^{(n)}$, factors
through $\T_X^{(n)}$.
\end{list}
\label{thm:main:higher}
\end{Thm}
We may add to (a) that in terms of the covering $\{V,X-p\}$ of $X$
as above, $\ell_1^{(n)}((Z_1\ldots Z_n)|_0)$ can be also described as
follows: it is a map
$$
H^0(X,\omega_X)\longrightarrow H^1(X,{\cal O}_X)
$$
sending the class represented by a form $\omega$ on $V$ to the class
represented by the function
$$
(-1)^n\zeta_n\contr\pounds_{\zeta_{n-1}}\ldots \pounds_{\zeta_1}
\omega
$$
on $V-p$. Here we assume that the lifts $\zeta_i\in
\d$ of $Z_i\in\Gamma(S,\Theta_S)$ converge and define regular vector
fields on
$V-p$.

\

\

\noindent Department of Mathematics\\Columbia University\\
New York, NY 10027

\

\

\noindent {\em Current address}:

\

\noindent Department of Mathematics\\Yeshiva University\\
500 W 185 St.\\
New York, NY 10033\\ \ \\{\em yk@yu1.yu.edu}

\end{document}